\documentclass[11pt]{amsart}

\newtheorem{theorem}{Theorem}[section]
\newtheorem{lem}[theorem]{Lemma}
\newtheorem{prop}[theorem]{Proposition}

\theoremstyle{definition}

\theoremstyle{remark}
\newtheorem{remark}[theorem]{Remark}
\numberwithin{equation}{section}

\newcommand{\half}{(1/2)}

\newcommand{\ra}{\rightarrow}

\newcommand{\C}{{\mathbb{C}}}
\newcommand{\HH}{{\mathbb{H}}}
\newcommand{\CP}{{\mathbb{CP}}}
\newcommand{\CaP}{{\mathbb{C}a\mathbb{P}}}
\newcommand{\HP}{{\mathbb{HP}}}
\newcommand{\RP}{{\mathbb{RP}}}

\newcommand{\KK}{\mathbb{K}}
\newcommand{\KP}{\mathbb{KP}}
\newcommand{\Z}{\mathbb{Z}}

\newcommand{\dds}{\frac{\text{d}}{\text{ds}}}

\newcommand{\lb}{\langle}
\newcommand{\rb}{\rangle}

\newcommand{\mg}{\mathfrak{g}}
\newcommand{\mh}{\mathfrak{h}}
\newcommand{\mk}{\mathfrak{k}}
\newcommand{\mm}{\mathfrak{m}}
\renewcommand{\mp}{\mathfrak{p}}
\newcommand{\bs}{\backslash}

\newcommand{\Add}{\text{Ad}}

\newcommand{\R}{\mathbb{R}}

\newcommand{\spp}{\text{sp}}

\begin{document}

\newcommand{\spacing}[1]{\renewcommand{\baselinestretch}{#1}\large\normalsize}
\spacing{1.14}

\title{Quasi-positive curvature on homogeneous bundles}

\author {Kristopher Tapp}

\address{Department of Mathematics\\ Bryn Mawr College\\
Philadelphia, PA 19010}
\email{ktapp@brynmawr.edu}



\begin{abstract}
We provide new examples of manifolds which admit a Riemannian metric with sectional curvature nonnegative, and strictly positive at one point.  Our examples include the unit tangent bundles of $\CP^n$, $\HP^n$ and $\CaP^2$, and a family of lens space bundles over $\CP^n$.  All new examples are consequences of a general sufficient condition for a homogeneous fiber bundle over a homogeneous space to admit such a metric.
\end{abstract}

\maketitle


\section{Introduction}\label{intro}
There are very few known examples of compact manifolds with strictly positive sectional curvature.  However, new examples have been recently constructed of nonnegatively curved manifolds with positive curvature either at a point (called \textit{quasi-positive curvature}) or on an open dense set of points (called \textit{almost-positive curvature}).  Gromoll and Meyer discovered a 7-dimensional exotic sphere with quasi-positive curvature~\cite{GM}.  This exotic sphere was later shown to admit almost-positive curvature~\cite{Wil},\cite{E2}.  Petersen and Wilhelm endowed $T^1S^4$ and a 6 dimensional quotient of $T^1S^4$ with almost-positive curvature~\cite{PW}.

More recently, in~\cite{Wilking}, Wilking discovered several families of manifolds admitting almost-positive curvature, including the projective tangent bundles $P_{\R}T\RP^n$, $P_{\C}T\CP^n$ and $P_{\HH}T\HP^n$, and a family of lens space bundles over $\CP^n$.  His result for $P_{\R}T\RP^n$ implies that its cover, $T^1S^n$, admits almost-positive curvature, which is particularly interesting in the cases of $T^1S^3=S^3\times S^2$ and $T^1 S^7=S^7\times S^6$.  Amongst his examples are spaces known not to admit positive curvature.  On the other hand, it remains unknown whether every manifold admitting quasi-positive curvature must admit almost-positive curvature. 

The main results of this paper are the following new examples:
\begin{theorem}\label{T:Ex} The following manifolds admit metrics with quasi-positive curvature:
\begin{enumerate}
\item The unit tangent bundles $T^1\CP^n$, $T^1\HP^n$ and $T^1\CaP^2$.
\item The homogeneous space $M=U(n+1)/\{\text{diag}(z^k,z^l,A)\mid z\in U(1),A\in U(n-1)\}$, where $(k,l)$ is a pair of integers with $k\neq 0$, and $n\geq 2$.
\item The homogeneous space $Sp(n+1)/\{\text{diag}(z,1,A)\mid z\in Sp(1), A\in Sp(n-1)\}$, $n\geq 2$.
\end{enumerate}
\end{theorem}

The space in part (2) is a lens space bundles over $\CP^n$.  Wilking proved that the sub-family with $k\cdot l<0$ admit almost-positive curvature.  Our larger family contains new examples, including the case $k=l=1$, which is $T^1\CP^n$.

The space in part (3) is an $S^{4n-1}$-bundle over $\HP^n$.  Wilking proved that a quotient of this space by a free $S^3$-action, namely the bi-quotient $Sp(1)\bs Sp(n+1)/Sp(1)\cdot Sp(n-1)$, admits almost-positive curvature.  Since our metric is $S^3$-invariant, it at least follows from our result that this quotient admits quasi-positive curvature.  More generally, our methods provide a simple way to prove quasi-versions of all of Wilking's almost-positive curvature results.

All of our examples are homogeneous bundles over homogeneous spaces.  If $B=H\bs G$ is a homogeneous space, then a fiber bundle $F\hookrightarrow M\ra B$ is called \emph{homogeneous} if the transitive right $G$-action on $B$ lifts to $M$.  Said differently, a homogenous fiber bundle over $B$ is one which can be written as a quotient $M=H\bs(G\times F)$, for some left action of $H$ on the fiber $F$.  If the action of $H$ on $F$ is transitive, then $M$ is diffeomorphic to $H\bs (G\times (H/K)) = G/K$, where $K$ is the isotropy group of this action.  So the bundle looks like:
$$H/K\ra G/K\ra G/H.$$  In this case, we will endow $M=H\bs(G\times F)$ with the submersion metric (using a natural left-invariant metric on $G$ and the product metric on $G\times F$), and derive conditions under which this metric on $M$ has quasi-positive curvature.

The author is pleased to thank Burkhard Wilking and Wolfgang Ziller for helpful conversations about this work.


\section{Summary of Conditions}
In this section, we summarize our conditions under which a homogeneous bundle admits quasi-positive curvature.
We adopt the following notation and assumptions for the remainder of the paper.  Let $B=H\bs G$ denote a homogenous space, with $G$ and $H$ compact Lie groups.  Let $g_0$ denote a bi-invariant metric on $G$, and assume that $H\bs(G,g_0)$ has positive curvature.  Let $F$ denote a compact Riemannian manfold which has positive curvature or is one dimensional.  Assume that $H$ acts transitively and isometrically on $F$ on the left.  Let $K\subset H$ denote the isotropy group at some point $p_0\in F$.  Denote the Lie algebras of $K\subset H\subset G$ as $\mk\subset\mh\subset\mg$.  Let $\mm:=\mh\ominus\mk$ and $\mp:=\mg\ominus\mh$, where ``$\ominus$'' denotes the $g_0$-orthogonal compliment.

Let $g_l^H$ denote the left-invariant and right-H-invariant metric on $G$ obtained from $g_0$ by rescalling in the direction of $\mh$.  More precisely, fix $t\in(0,1)$ and define:
\begin{equation}\label{gl}
g_l^H(X,Y)=g_0(X^{\mp},Y^{\mp})+t\cdot g_0(X^{\mh},Y^{\mh})
\end{equation}
where $X^{\mp}$ (respectively $X^{\mh}$) denotes the projection of $X$ orthogonal to (respectively onto) $\mh$.  Notice that $(G,g_l^H)$ is nonnegatively curved by~\cite{Eschenburg}, since it can be described as a submersion metric:
\begin{equation}\label{lambda}(G,g_l^H) = ((G,g_0)\times(H,\lambda\cdot g_0|_H))/H,\end{equation}
where $\lambda = t/(1-t)$.  

Let $M = H\bs((G,g_l^H)\times F)$, where $H$ acts diagonally, and denote the projections as follows:
$$(G,g_l^H)\times F\stackrel{\pi}{\ra} M \stackrel{\phi}{\ra} H\bs (G,g_l^H).$$
Notice that the maps $\pi$ and $\phi\circ\pi$ are Riemannian submersions by construction, so it follows that $\phi$ is a Riemannian submersion as well.  The fibers of $\phi$ are not in general totally geodesic, although they would be if $g_l^H$ were replaced by a right-invariant and left-$H$-invariant metric. 

As for the isometries remaining on $M$, the following is straightforward to verify:
\begin{remark}\label{REM} The right $H$-action on the first factor of $(G,g_l^H)\times F$ induces an isometric action of $H$ on $M$.  A subgroup $L\subset H$ acts freely if and only if $L\cap(g^{-1}\cdot K\cdot g)=\{0\}$ for all $g\in G$.  In this case, the quotient $M/L$ is diffeomorphic to the bi-quotient $L\bs G/K$.
\end{remark}

The space $M$ has nonnegative curvature.  Our first goal is to derive conditions under which $M$ has points of positive curvature.  Our conditions turn out not to depend on $t$.
\begin{theorem}\label{L:C2}\hspace{.1in}
\begin{enumerate}
\item If $[Z,W]\neq 0$ for all linearly independent vectors $Z\in\mg\ominus\mk$ and $W\in\mp$, then $M$ has positive curvature.
\item
If there exists $A\in\mg$ such that $[Z^{\mh},[A,W]^{\mh}]\neq 0$ for all linearly independent vectors $Z\in\mg\ominus\mk$ and $W\in\mp$ for which $[Z,W]=0$, then $M$ has quasi-positive curvature.
\item
If $(G,H)$ is a compact rank one symmetric pair, and if there exists $A\in\mp$ such that $[X,A]\neq 0$ for all non-zero $X\in\mm$, then $M$ has quasi-positive curvature.
\end{enumerate}
\end{theorem}
We will prove this theorem in section~\ref{S:proof}.  We prove part (1) by deriving straightforward conditions for zero-curvature planes.  We prove part (2) by differentiating these conditions to show that for small $\epsilon>0$, points in $M$ of the form $\pi((\exp(-\epsilon\cdot A),p_0))$ have positive curvature.  When $(G,H)$ is a symmetric pair, the hypothesis of part (2) simplifies to the hypothesis of part (3).  All of our new examples are consequences of part (3) of the theorem.  We construct these new examples in section~\ref{S:E}, and also recover quasi-versions of Wilking's almost-positive curvature results.  In fact, our metrics are isometric to his, which we prove in section~\ref{S:biquotient}.

\section{Positive Curvature?}

At first glance, it seems possible to use part 1 of Theorem~\ref{L:C2} to find new examples of positively curved manifolds.  However, the hypothesis implies that $[Z,W]\neq 0$ for all non-zero vectors $Z\in\mm$ and $W\in\mp$, which is called ``fatness'' of the homogeneous bundle (the hypothesis if equivalent to fatness when $(G,H)$ is a rank one symmetric pair).  See~\cite{fat} for an overview of literature related to the fatness condition.

Berard Bergery classified all fat homogeneous bundles in~\cite{BB}.  Our theorem provides hindsight motivation for his classification.  If $(G,H)$ is a rank one symmetric pair, he found that the bundle is fat if and only if $M$ admits a homogeneous metric of positive curvature.  Further, if the fiber dimension is $>1$, he found that fatness implies that $(G,H)$ must be a rank one symmetric pair.  Therefore, no new examples of positive curvature can be found with Theorem~\ref{L:C2}, with the possible exception of circle bundles.

It is already known which circle bundles over rank one symmetric spaces admit positive curvature.  The three non-symmetric positively curved normally homogeneous spaces discovered in~\cite{Berger} and~\cite{WB} are odd dimensional, making a circle bundle over one be even-dimensional.  But our metric on a circle bundle looks like $M=H\bs ((G,g_l^H)\times S^1)$, which admits the free isometric $S^1$-action induced by the action of $S^1$ on the second factor of $(G,g_l^H)\times S^1$.  So positive curvature on $M$ would contradict Berger's Theorem, which says that an even dimensional positively curved manifold cannot admit a nonvanishing Killing field.  Thus, our construction does not yield new examples of positive curvature.


\section{Proof of Theorem~\ref{L:C2}}\label{S:proof}
Notice that any point of $M$ can be written as $\pi((g^{-1},p_0))$ for some $g\in G$.  If $(g^{-1},p_0)$ does not contain a $\pi$-horizontal zero-curvature plane, then $\pi((g^{-1},p_0))$ is a point of positive curvature.  We begin with the following lemma, which is more generally valid when $g_l^H$ is replaced by any left-invariant metric on $G$.

\begin{lem}\label{L1}
There exists a $\pi$-horizontal zero-curvature plane at $(g^{-1},p_0)$ if and only if there exist vectors
\begin{gather*}
X\in \mg\ominus\Add_g(\mk)=\{Y\in\mg\mid g_l^H(Y,\Add_gA)=0\,\,\,\forall A\in \mk\},\\
W_i\in\mg\ominus \Add_g(\mh)=\{Y\in\mg\mid g_l^H(Y,\Add_gA)=0\,\,\,\forall A\in \mh\}
\end{gather*}
such that $\text{span}\{X+W_1,X+W_2\}$ is a 2-plane in $\mg$ with zero-curvature with respect to $g_l^H$.
\end{lem}
\begin{proof}

Suppose that $\sigma$ is a $\pi$-horizontal zero-curvature plane at $(g^{-1},p_0)$.  Since $\sigma$ has zero-curvature,
it is spanned by vectors of the form $\{dL_{g^{-1}}\overline{W}_1+V,dL_{g^{-1}}\overline{W}_2+V\}$, where 
$V\in T_{p_0} F$, and $\overline{W}_1,\overline{W}_2\in\mg$ span a $g_l^H$-zero curvature plane (notice that $\overline{W}_1$ and $\overline{W}_2$ are linearly independent because $V$ is not $\pi$-horizontal).

The vertical space of $\pi$ at $(g^{-1},p_0)$ is the set $\{dR_{g^{-1}} A + A(p_0)\mid A\in\mh\}$,
where $A(p_0)\in T_{p_0} F$ denotes the value at $p_0$ of the Killing field on $F$ associated with $A$.  So for all $i\in\{1,2\}$ and $A\in\mh$, we have:
$$0=\lb dL_{g^{-1}}\overline{W}_i+V,dR_{g^{-1}} A + A(p_0)\rb = g_l^H( \overline{W}_i,Ad_g A ) + \lb V,A(p_0)\rb.$$
In particular,$\overline{W}_1,\overline{W}_2\in\mg\ominus\Add_g(\mk)$.  Further, $\overline{W}_1$ and $\overline{W}_2$ have the same $g_l^H$-orthogonal projections onto $\Add_g(\mh)$.  We denote their common projection
as $X$.  Denote $W_i=\overline{W}_i-X$.  Then $\text{span}\{X+W_1,X+W_2\}=\text{span}\{\overline{W}_1,\overline{W}_2\}$
is a zero-curvature plane in $\mg$ as required.

Conversely, suppose there exist vectors $X\in\mg\ominus\Add_g(\mk)$ and $W_i\in\mg\ominus\Add_g(\mh)$ such that
$\text{span}\{X+W_1,X+W_2\}$ is a zero-curvature plane in $\mg$.  For any $V\in T_{p_0} F$, the plane $\sigma=\text{span}\{dL_{g^{-1}}(X+W_1) + V,dL_{g^{-1}}(X+W_2) + V\}$ is a zero-curvature plan in $G\times F$ at $(g^{-1},p_0)$.  It will suffice to choose $V$ such that $\sigma$ is $\pi$-horizontal; that is, such that for all $i\in\{1,2\}$ and $A\in\mh$,
$$0=\lb dL_{g^{-1}}(X+W_i) + V, dR_{g^{-1}} A + A(p_0)\rb =  g_l^H(X,\Add_gA) + \lb V,A(p_0)\rb.$$
By linearity, it will suffice to find $V\in T_{p_0}F$ such that $\lb V,A_i(p_0)\rb = -g_l^H(X,\Add_gA_i)$
for each element $A_i$ of a basis of $\mh$.  Choose a basis $\{A_i\}$ such that $A_i\in\mk$ for $i\leq\text{dim}(\mk)$ and such that $\{A_i(p_0)\mid i>\text{dim}(\mk)\}$ is a basis of $T_{p_0}F$.  Since $X$ is orthogonal to $\mk$, it is easy to see that such a vector $V$ can be chosen.
\end{proof}

In order to apply Lemma~\ref{L1}, we require Eschenburg's description of the planes in $\mg$ which have zero curvature with respect to $g_l^H$.  The case $t=1/2$ is found in~\cite{Eschenburg}, and the general case is similar.  To phrase his condition, we define $\Phi:\mg\ra\mg$ by the rule: $\Phi(A)=t\cdot A^{\mh}+A^{\mp}.$
It's easy to verify that for all $A,B\in\mg$,
$$g_l^H(A,B)=g_0(A,\Phi(B))\text{\quad and \quad} g_0(A,B) = g_l^H(A,\Phi^{-1}(B)).$$

\begin{lem}[Eschenburg]\label{esc} Let $\sigma=\text{span}\{X,Y\}$ be a plane in $\mg$.
\begin{enumerate}
\item $\sigma$ has zero-curvature if and only if $[\Phi(X),\Phi(Y)]=0$ and $[X^{\mh},Y^{\mh}]=0$.
\item If $(G,H)$ is a symmetric pair, then $\sigma$ has zero-curvature if and only if $[X,Y]=0$ and $[X^{\mh},Y^{\mh}]=0$.
\end{enumerate}
\end{lem}

Combining Lemmas~\ref{L1} and~\ref{esc} yields the following condition:

\begin{lem}\label{L:C1}
There exists a $\pi$-horizontal zero-curvature plane at $(g^{-1},p_0)$ if and only if there exist linearly independent vectors $Z\in\mg\ominus\mk$ and $W\in\mp$ such that $[Z,W]=0$ and $[(\Add_gZ)^{\mh},(\Add_g W)^{\mh}]=0$.
\end{lem}

\begin{proof}

Assume there exists a $\pi$-horizontal zero curvature plane at $(g^{-1},p_0)$. By Lemma~\ref{L1},
there exists vectors $\overline{X}\in\mg\ominus\Add_g(\mk)$ and $\overline{W_i}\in\mg\ominus\Add_g(\mh)$ such that $$\text{span}\{\overline{X}+\overline{W}_1,\overline{X}+\overline{W}_2\} =\text{span}\{\overline{X}+\overline{W}_1,\overline{W}_2-\overline{W}_1\}$$
is a zero curvature plane in $\mg$ with respect to $g_l^H$.

It is possible to write $\overline{X}+\overline{W}_1=\Phi^{-1}(\Add_gZ)$ and $\overline{W}_2-\overline{W}_1=\Phi^{-1}(\Add_g W)$ for some $Z,W\in\mg$.  In fact, $Z\in\mg\ominus\mk$ because for all $A\in\mk$,
$$0=g_l^H(\overline{X}+\overline{W}_1,\Add_gA) = g_l^H(\Phi^{-1}(\Add_gZ),\Add_gA)=g_0(\Add_gZ,\Add_gA)=g_0(Z,A).$$
Similarly, $W\in\mp$.  Applying Lemma~\ref{esc} gives:
\begin{eqnarray*}
0 & = & [\Phi(\overline{X}+\overline{W}_1),\Phi(\overline{W}_2-\overline{W}_1)]
   =  [\Add_g(Z),\Add_g(W)] = [Z,W].\\
0 & = & [(\overline{X}+\overline{W}_1)^{\mh},(\overline{W}_2-\overline{W}_1)^{\mh}]
       =[(1/t)(\Add_g(Z))^{\mh},(1/t)(\Add_g W)^{\mh}].
\end{eqnarray*}
Therefore, the vectors $\{Z,W\}$ satisfy the conclusions of the lemma.  The other direction of the lemma follows analagously.
\end{proof}

The previous lemma simplifies with the symmetric pair conditions $[\mp,\mp]\subset\mh$ and $[\mp,\mh]\subset\mp$.

\begin{lem}\label{mainL}
Assume that $(G,H)$ is a symmetric pair.  Then there exists a $\pi$-horizontal zero-curvature plane at $(g^{-1},p_0)$ if and only if there exists non-zero vectors $X\in\mm$ and $W\in\mp$ such that $[X,W]=0$ and $[(\Add_gX)^{\mh},(\Add_gW)^{\mh}]=0$.
\end{lem}

\begin{proof}
Assume there exists a $\pi$-horizontal zero curvature plane at $(g^{-1},p_0)$.  So by Lemma~\ref{L:C1},
there exists linearly independent vectors $Z\in\mg\ominus\mk$ and $W\in\mp$ for which $[Z,W]=0$ and $[(\Add_gZ)^{\mh},(\Add_g W)^{\mh}]=0$. Then, $0=[Z^{\mh}+Z^{\mp},W]=[Z^{\mh},W]+[Z^{\mp},W]$.  But since $[Z^{\mh},W]\in\mp$ and $[Z^{\mp},W]\in\mh$, both vectors must be zero.  Since $H\bs(G,g_o)$ has positive curvature, $Z^{\mp}$ and $W$ are parallel.  So $Z^{\mh}\neq 0$, and the pair $\{X=Z^{\mh},W\}$ satisfies the conclusions of Lemma~\ref{mainL}.  The other direction of the lemma is argued similarly.
\end{proof}

Finally, we prove Theorem~\ref{L:C2}.
\begin{proof}[Proof of Theorem~\ref{L:C2}]
Part (1) is immediate from Lemma~\ref{L:C1}.  To prove part (2), assume that $A\in\mg$ satisfies its hypothesis.  Notice that Lemma~\ref{L:C1} remains true if the phrase ``linearly independent'' is replaced by the phrase ``$g_0$-orthonormal''.  So fix $g_0$-orthonormal vectors $Z\in\mg\ominus\mk$ and $W\in\mp$ for which $[Z,W]=0$.  Define:
$$f(s)=\left|\left[\left(\Add_{\exp(sA)} Z\right)^{\mh},\left(\Add_{\exp(sA)}W\right)^{\mh}\right]\right|^2.$$
Then $f(0)=0$, $f'(0)=0$ and:
\begin{eqnarray*}
\frac{1}{2}f''(0)
 & = & \left|\dds\Big|_{s=0}\left[\left(\Add_{\exp(sA)}Z\right)^{\mh},\left(\Add_{\exp(sA)}W\right)^{\mh}\right]\right|^2\\
 & = & \left|\left[Z^{\mh},\dds\Big|_{s=0}\left(\Add_{\exp(sA)}W\right)^{\mh}\right]\right|^2\
  =  \left|\left[Z^{\mh},\left(\dds\Big|_{s=0}\Add_{\exp(sA)}W\right)^{\mh}\right]\right|^2\\
 & = & \left|\left[Z^{\mh},\left[A,W\right]^{\mh}\right]\right|^2>\delta>0.
\end{eqnarray*}
By compactness (of the space of orthonormal vectors $\{Z,W\}$ with $[Z,W]=0$), $\delta$ can be chosen to depend only on $A$ (not on $Z$ and $W$).  It follows that $\delta'>0$ can be chosen (independent of $Z$ and $W$) such that $f(s)>0$ for all $s\in(0,\delta']$.  It now follows from Lemma~\ref{L:C1} that for $s\in(0,\delta']$, $\pi((\exp(-s\cdot A),p_0))$ is a point of positive curvature.

To prove part (3), suppose there exists $A\in\mp$ such that $[X,A]\neq 0$ for all non-zero $X\in\mm$.  Then for all non-zero $X\in\mm$ we have: $$0\neq g_0([X,A],[X,A]) = g_0([X,[X,A]],A).$$ In particular, $[X,[X,A]]\neq 0$.

Now let $Z\in\mg\ominus\mk$ and $W\in\mp$ be linearly independent vectors for which $[Z,W]=0$.  Using the assumption that $(G,H)$ is a rank one symmetric pair, $X:=Z^{\mh}\neq 0$ and $[X,W]=0$.  From part (2) of the theorem, it will suffice to show that the following is non-zero:
$$[Z^{\mh},[A,W]^{\mh}] = [X,[A,W]] = -[W,[X,A]] \text{ (Jacobi identity)}$$
But if $[W,[X,A]]=0$, then $W$ would be parallel to $[X,A]$, which would mean $[X,[X,A]]=0$.  This is a contradiction.  Therefore $M$ has quasi-positive curvature.
\end{proof}
\section{The Space $S^2\times S^3$}
In this section, we use Lemma~\ref{mainL} to recover Wilking's theorem that $S^2\times S^3$ admits a metric with almost-positive curvature.  Although our proof is very similar to the original, we find it worthwhile to translate the original proof into the vocabulary of this paper.

\begin{prop}[Wilking]\label{S2S3} The space $S^2\times S^3$ admits almost-positive curvature.
\end{prop}
\begin{proof}
One way to describe $S^3$ as a symmetric space $S^3=H\bs G$ is $G=S^3\times S^3$ and $H=\Delta(S^3)$.  If we let $H$ act on $(S^2,\text{round})$ as $SO(3)$, then the associated bundle $M=H\bs(G\times S^2)$ is diffeomorphic to $T^1S^3=S^3\times S^2$.  One isotropy group of the action of $H$ on $S^2$ is 
$K=\{(e^{i\theta},e^{i\theta})\mid \theta\in S^1\}$.  The Lie algebras of $K\subset H\subset G$ are: $$\mk=\text{span}\{(i,i)\}\subset\mh=\Delta(\spp(1))\subset\mg=\spp(1)\times\spp(1).$$

Let $g=(g_1,g_2)\in G$. Let $X=(a,a)\in\mm=\mh\ominus\mk$ and $W=(b,-b)\in\mp=\mg\ominus\mh$ be non-zero vectors, which means that $a,b\in\spp(1)$ are not zero, and $a$ is orthogonal to $i$.  Assume that $[X,W]=0$ and $[(\Add_gX)^{\mh},(\Add_gW)^{\mh}]=0$.  By Lemma~\ref{mainL}, it will suffice to prove that $g$ lies in the compliment of an open dense set of $G$.

Since $0=[X,W]=[(a,a),(b,-b)]=([a,b],[a,-b])$, we see that $[a,b]=0$.  Thus $a$ is parallel to $b$, so $b=\lambda a$ for some $\lambda\in\R$.  Next,
\begin{eqnarray*}
0 & = & [(\Add_g(a,a))^{\mh},(\Add_g(\lambda a,-\lambda a))^{\mh}] \\
  & = & \lambda^2[(\Add_{g_1}a,\Add_{g_2}a)^{\mh},(\Add_{g_1}a,-\Add_{g_2}a)^{\mh}]\\
  & = & \lambda^2\left[\left(\frac{\Add_{g_1}a + \Add_{g_2}a}{2},\frac{\Add_{g_1}a + \Add_{g_2}a}{2}\right),\left(\frac{\Add_{g_1}a - \Add_{g_2}a}{2},\frac{\Add_{g_1}a - \Add_{g_2}a}{2}\right)\right]
\end{eqnarray*}
From this we see that $[\Add_{g_1}a,\Add_{g_2}a]=0$.  This implies that $\Add_{g_1}a = \pm \Add_{g_2}a$; in other words, $a=\pm\Add_{(g_1^{-1}g_2)}a$.  Using the fact that $a\perp i$, this implies that $(g_1^{-1}g_2)\perp 1$ or $(g_1^{-1}g_2)\perp i$.
\end{proof}

\section{Old and New examples}\label{S:E}
In this section, we will use part (3) of Theorem~\ref{L:C2} to obtain all of the new examples in Theorem~\ref{T:Ex}, and also to recover quasi-versions of Wilking's almost-positive curvature results.

\begin{prop} The unit tangent bundle of any compact rank one symmetric space admits quasi-positive curvature.
\end{prop}
\begin{proof} Assume $(G,H)$ is a compact rank one symmetric pair.  Let $S=G/H$, which is one of $\RP^n, S^n, \HP^n$ or $\CaP^2$.  There is a natural $G$ action on the unit tangent bundle $T^1 S$, obtained by differentiating the $G$-action on $S$.  This action on $T^1 S$ is transitive because the isotropy representation of $H$ on $\mp=T_{[e]}S$ is transitive on the unit-sphere $\mp^1$, which can be checked for each of the four possibilities for $S$.  Fix $A\in\mp^1$.  Then $T^1 S=G/K$, where $K$ is the collection of elements of $H$ fixing $A$ under the isotropy representation of $H$ on $\mp$.  By part (3) of Theorem~\ref{L:C2}, it suffices to prove that $[A,X]\neq 0$ for all non-zero $X\in\mm$.

Each element, $Y\in\mh$ determines a Killing vector field on $\mp^1$ via the isotropy representation of $H$ on $\mp$.  The value of this Killing field at $V\in\mp^1$ equals $[Y,V]$.  Since the action of $H$ on $\mp^1$ is transitive, these Killing fields must span $T_V\mp^1$ for every $V\in\mp^1$.  Since the Killing field associated to an element $Y\in\mk$ vanishes at $A$, and since $\text{dim}(\mm)=\text{dim}(\mp^1)$, the Killing fields associated to a basis for $\mm$ must be linearly independent at $A$.  Therefore, $[A,X]\neq 0$ for all non-zero $X\in\mm$.
\end{proof}

\begin{prop} When $n$ is even, there is a free $S^1$-action on $T^1 S^n$ such that the quotient, $T^1S^n/S^1 = SO(2)\bs SO(n+1)/SO(n-1)$, admits quasi-positive curvature.
\end{prop}
Wilking proved that $T^1S^n/S^1$ admits almost-positive curvature.
\begin{proof}
The groups for $T^1S^n$ are $K=SO(n-1)\subset H=SO(n)\subset G=SO(n+1)$.  By the previous proposition, $T^1 S^n$ admits quasi-positive curvature.  By Remark~\ref{REM}, $H=SO(n)$ acts isometrically on $T^1 S^n$, and the subgroup $L=SO(2)$ embedded diagonally in $H$ acts freely.  Therefore the quotient admits quasi-positive curvature.
\end{proof}

\begin{prop} The projective tangent bundles $P_{\R}T\RP^n$, $P_{\C}T\CP^n$ and $P_{\HH}T\HP^n$ admit metrics with quasi-positive curvature.
\end{prop}

Wilking proved that these spaces admit almost-positive curvature.  The projective tangent bundle of $\CaP^2$ is known to admit a homogeneous metric with positive curvature~\cite{Wallach}.

\begin{proof}
Let $\KK\in\{\R,\C,\HH\}$.  Let $G(n)$ denote $O(n)$, $U(n)$ or $Sp(n)$, depending on $\KK$.  The unit tangent bundle $T^1\KP^n$ and the projective tangent bundle $P_{\KK}T\KP^n$ come respectively from the groups:
\begin{gather*}
\{\text{diag}(z,z,A)\mid z\in G(1), A\in G(n-1)\}\subset G(1)\cdot G(n)\subset G(n+1) \\
\{\text{diag}(z_1,z_2,A)\mid z_i\in G(1), A\in G(n-1)\}\subset G(1)\cdot G(n)\subset G(n+1).
\end{gather*}

In the previous proposition, we verified that the groups $K\subset H\subset G$ for $T^1\KP^n$ satisfy the condition for quasi-positive curvature.  Since $K$ is strictly larger in the groups for $P_{\KK}T\KP^n$, the condition is also satisfied there.
\end{proof}

\begin{prop}\label{PPP}
The homogeneous space $M_{kl}=U(n+1)/K_{kl}$ admits quasi-positive curvature, where $(k,l)$ is a pair of integers with $k\neq 0$, $n\geq 2$ and $K_{kl}=\{\text{diag}(z^k,z^l,A)\mid z\in S^1,A\in U(n-1)\}$.
\end{prop}

When $l=0$, we lose no generality in assuming that $k=1$.  When $l\neq 0$, we lose no generality in assuming that $k$ and $l$ are relatively prime, since dividing both by a common factor does not change the group $K_{kl}$.  In these case, $M$ is a homogeneous bundle over $\CP^n=U(n+1)/(U(1)\cdot U(n))$ with fiber equal to the homogeneous lens space $L_k=S^{2n-1}/\Z^k=U(n)/(\Z_k\cdot U(n-1))$.  To describe the transitive isometric action of $H=U(1)\cdot U(n)$ on $L_k$ which yields these bundles, denote a point of  $L_k$ as $[C]$, where $C\in U(n)$.  The action is:
$$[C]\stackrel{(z,A)}{\mapsto} [z^{-l/k}\cdot A\cdot C],$$
where $z^{-l/k}$ means the $-l^{\text{th}}$ power of any $k^{\text{th}}$ root of $z$.  Since the answer does not depend on the choice of root, this action is well-defined, and has isotropy group $K_{kl}$.

A more general action of $H$ on $L_k$ is: $$[C]\stackrel{(z,A)}{\mapsto} [z^{-l_1/k}\cdot(\det A)^{l_2/k}\cdot A\cdot C].$$  However, the total space of the resultant lens space bundle over $\CP^n$ depends only on the integers $k$ and $l_1+l_2$, so we lose no generality in assuming that $l_2=0$.  So see this, notice that isotropy group of this more general action is:
\begin{eqnarray*}
K_{k l_1 l_2} & = & \{(z,w,A)\in U(1)\cdot U(1)\cdot U(n-1)\mid z^{-l_1/k}\cdot w^{l_2/k}\cdot (\det A)^{l_2/k}\cdot w\in \Z_k\} \\
  & = & \{(z,w,A)\mid z^{-l_1}\cdot w^{k+l_2}\cdot(\det A)^{l_2}=1\}.
\end{eqnarray*}
Next observe that $SU(n+1)\subset U(n+1)$ acts transitively on $M_{k l_1 l_2}=U(n+1)/K_{k l_1 l_2}$ because every coset intersects $SU(n+1)$.  So $M_{k l_1 l_2}$ is diffeomorphic to $SU(n+1)/K'_{k l_1 l_2}$, where:
$$K'_{k l_1 l_2}=SU(n+1)\cap K_{k l_1 l_2} = \{(z,w,A)\in SU(n+1)\mid z^{-l_1-l_2}\cdot w^k = 1\}.$$

\begin{proof}[Proof of Proposition~\ref{PPP}]
The Lie algebra of $H$ is $\mh=\mathfrak{u}(1)\cdot\mathfrak{u}(n)$.  The Lie algebra of $K$ is:
$$\mk = \{\text{diag}(tki,tli,B)\in\mathfrak{u}(1)\cdot\mathfrak{u}(1)\cdot\mathfrak{u}(n-1)\mid t\in\R\}.$$
We use the bi-invariant metric on $\mg=\mathfrak{u}(n+1)$ defined by $\lb A,B\rb = \text{Real}(\text{trace}(A\cdot\overline{B}^{T}))$.  Arbitrary vectors $X\in\mm$ and $A\in\mp$ look like:
$$
X=\left(\begin{matrix} -tli & 0 & 0 & \cdots & 0 \\ 0 & tki & Z_1 & \cdots & Z_{n-1} \\
                       0 & -\overline{Z}_1 & 0 & \cdots & 0 \\ \vdots & \vdots & \vdots & & \vdots \\
                       0 & -\overline{Z}_{n-1} & 0 & \cdots & 0 \end{matrix}\right), \qquad
A=\left(\begin{matrix} 0 & W_1 & \cdots & W_n \\ -\overline{W}_1 & 0 & \cdots & 0 \\ \vdots & \vdots & & \vdots \\ -\overline{W}_n & 0 & \cdots & 0 \end{matrix}\right),$$
where $t\in\R$ and $W_i,Z_i\in\C$.  Using the shorthand $A=\{W_1,W_2,...,W_n\}$,
$$[X,A]=\{-t(l+k)iW_1 + \overline{Z}_1W_2+\cdots+\overline{Z}_{n-1}W_n,-Z_1W_1-tliW_2,...,-Z_{n-1}W_1-tliW_n\}$$
If we set $W_1=W_2=1, W_3=\cdots=W_n=0$, then $[X,A]=0$ implies $X=0$.  Thus, no non-zero vector in $\mm$ commutes with $A$, so $M_{kl}$ admits quasi-positive curvature.

\end{proof}

\begin{prop} The homogeneous space $M=Sp(n+1)/\{\text{diag}(z,1,A)\mid z\in Sp(1), A\in Sp(n-1)\}$ admits quasi-positive curvature.  Further, there is a free isometric $S^3$-action $M$, so the quotient $M/S^3=Sp(1)\bs Sp(n+1)/Sp(1)\cdot Sp(n-1)$ admits quasi-positive curvature.
\end{prop}
Wilking proved that $M/S^3$ admits almost-positive curvature, but he did not prove anything about $M$ itself.  Notice that $M$ is the total space of a homogeneous $S^{4n-1}$-bundle over $\HP^n=Sp(n+1)/(Sp(1)\cdot Sp(n))$.  The corresponding action of $H=Sp(1)\cdot Sp(n)$ is the one whereby $p\in S^{4n-1}$ is sent by $(A,q)\in H$ to $A\cdot p$.  Notice that sending $p$ to $A\cdot p\cdot q^{-1}$ would yield $T^1\HP^n$.  

\begin{proof}
We use the bi-invariant metric on $\text{sp}(n+1)$ defined by $\lb A,B\rb = \text{Real}(\text{trace}(A\cdot\overline{B}^{T}))$.  Arbitrary elements $X\in\mm$ and $A\in\mp$ can be described as:
$$
X=\left(\begin{matrix} 0 & 0 & 0 & \cdots & 0 \\ 0 & Y & Z_1 & \cdots & Z_{n-1} \\
                       0 & -\overline{Z}_1 & 0 & \cdots & 0 \\ \vdots & \vdots & \vdots & & \vdots \\
                       0 & -\overline{Z}_{n-1} & 0 & \cdots & 0 \end{matrix}\right), \qquad
A=\left(\begin{matrix} 0 & W_1 & \cdots & W_n \\ -\overline{W}_1 & 0 & \cdots & 0 \\ \vdots & \vdots & & \vdots \\ -\overline{W}_n & 0 & \cdots & 0 \end{matrix}\right),$$
where $W_i,Z_i,Y\in\mathbb{H}$ with $\text{Real}(Y)=0$.  Using the shorthand $A=\{W_1,W_2,...,W_n\}$, we have:
$$[X,A]=\{-YW_1+W_2\overline{Z}_1 +\cdots+W_n\overline{Z}_{n-1},-W_1Z_1,...,-W_1Z_{n-1}  \}$$
The choice $W_1=1,W_2=\cdots=W_n=0$ insures that no non-zero vector in $\mm$ commutes with $A$.

By Remark~\ref{REM}, the subgroup $L=\{z\cdot I\mid z\in Sp(1)\}\subset H$ acts freely and isometrically on $M$, where $I$ denotes the identity in $Sp(n+1)$.  $L$ is isomorphic to $S^3$, and the quotient $M/S^3$ inherits quasi-positive curvature.
\end{proof}
\section{Normal biquotient metrics} \label{S:biquotient}
So far we have searched for points of positive curvature on the space $H\bs((G,g_l^H)\times F)$, which is diffeomorphic to the homogeneous space $G/K$.  Most of Wilking's examples are homogeneous spaces $G/K$ with bi-quotient metrics of the form:
\begin{equation}\label{local}\Delta(G)\bs((G,g_l^H)\times(G,g_l^H))/1\times K.\end{equation}
More precicely, his metrics on $P_{\R}T\RP^n$, $P_{\C}T\CP^n$, $P_{\HH}T\HP^n$, and $U(n+1)/K_{lk}$ all have this form.  Further, his example $SO(2)\bs SO(2n+1)/SO(2n-1)$ is a quotient of $T^1S^{2n}$ by a free isometric $S^1$-action, so his metric on this space also comes from one of the above form.  His only example not coming from a metric of the above form is the bi-quotient $Sp(1)\bs Sp(n+1)/Sp(1)\cdot Sp(n-1).$

In this section, we show that these biquotient metrics are isometric to our metrics.

\begin{prop}
The normal biquotient $\Delta(G)\bs((G,g_l^H)\times(G,g_l^H))/1\times K$ is isometric to the space $H\bs((G,\overline{g}_l^H)\times F)$, where $F$ has a normal homogeneous metric and $\overline{g}_l^H$ is defined like $g_l^H$, but with a different choice of $t$ and a different scaling of the bi-invariant metric $g_0$.

\end{prop}
In the proof we use the standard notational convention for biquotient.  That is, if $A\subset G\times G$, then $G//A$ denotes the orbit space of the action of $A$ on $G$ defined by $g\stackrel{(a_1,a_2)\in A}{\mapsto} a_1\cdot g\cdot a_2^{-1}$.  In the special case that $A=\{(a_1,a_2)\mid a_1\in A_1,a_2\in A_2\}$, we denote $G//A$ as $A_1\bs G/A_2$.
\begin{proof}
The equivalent definition of $g_l^H$ from equation~\ref{lambda} is:
$$(G,g_l^H) = ((G,g_0)\times(H,\lambda\cdot g_0))/H.$$
The isometry sends $[g,h]\in ((G,g_0)\times(H,g_0))/H$ to $g\cdot h^{-1}\in (G,g_l^H)$.  Notice that $(G,g_l^H)$ is left invariant and right H-invariant.  Via this isometry, the left G-action on $((G,g_0)\times(H,\lambda\cdot g_0))/H$ acts only on the first factor as $[g,h]\stackrel{a\in G}{\mapsto}[ag,h]$, and the right $H$-action acts only on the second factor as $[g,h]\stackrel{a\in H}{\mapsto}[g,a^{-1}h]$.  Using this, the biquotient 
$N =  \Delta(G)\bs((G,g_l^H)\times(G,g_l^H))/1\times K$ can be re-described as:
$$
\{(g,1,g,k^{-1})\mid g\in G,k\in K\}
  \bs((G,g_0)\times(H,\lambda\cdot g_0)\times(G,g_0)\times(H,\lambda\cdot g_0))/\{(h_1,h_1,h_2,h_2)\mid h_i\in H\}
$$

The two copies of $(G,g_0)$ in the above description of $N$ can be combined using the identity:
$$(G,\half g_0) = G\bs((G,g_0)\times(G,g_0)).$$
The isometry sends $[a,b]\in G\bs((G,g_0)\times(G,g_0))$ to $a^{-1}\cdot b\in (G,\half g_0)$.  Via this isometry, the right $G$-action on the first factor of $G\bs((G,g_0)\times(G,g_0))$ becomes the following right action of $G$ on $(G,\half g_0)$: $g\stackrel{a\in G}{\mapsto} a^{-1}g$.  Similarly, the right $G$-action on the second factor of $G\bs((G,g_0)\times(G,g_0))$ becomes the following right action of $G$ on $(G,\half g_0)$: $g\stackrel{a\in G}{\mapsto} ga$.  Therefore $N$ can be re-described as:
$$N=\{h_1^{-1},1,k^{-1}\}\bs((G,\half g_0)\times(H,\lambda\cdot g_0)\times(H,\lambda\cdot g_0))/\{h_2,h_1,h_2\},$$
which is our notational shorthand for the biquotient:
$$((G,\half g_0)\times(H,\lambda\cdot g_0)\times(H,\lambda\cdot g_0))//\{((h_1^{-1},1,k^{-1}),(h_2^{-1},h_1^{-1},h_2^{-1}))\mid k\in K,h_i\in H\}.$$

Next using the isometry of $(G,\half g_0)\times(H,\lambda\cdot g_0)\times(H,\lambda\cdot g_0)$ defined as $(a,b,c)\mapsto(a^{-1},b,c^{-1})$, we re-describe $N$ as:
$$N=\{h_2,1,h_2\}\bs((G,\half g_0)\times(H,\lambda\cdot g_0)\times(H,\lambda\cdot g_0))/\{h_1,h_1,k\}.$$

On the other hand, if $F$ has the normal homogeneous metric $F=(H,\lambda\cdot g_0)/K$, then:
$$H\bs((G,g_l^H)\times F) = \{h_2,1,h_2\}\bs((G,g_0)\times(H,\lambda\cdot g_0)\times(H,\lambda\cdot g_0))/\{h_1,h_1,k\}.$$
\end{proof}


\bibliographystyle{amsplain}

\begin{thebibliography}{9}

\bibitem{Berger} M. Berger, \emph{Les varietes Riemanniennes homogenes normales simplement connexes a Courbure strictment positive}, Ann. Scuola Norm. Sup. Pisa \textbf{30} (1975), 43-61.
\bibitem{BB} L. B\'{e}rard Bergery, \emph{Certaines fibrations d'espaces homog\`{e}nes Riemannienes}. Commposito Mathematica. \textbf{30} (1975), 43-61.
\bibitem{Eschenburg} J.-H.\ Eschenburg, \emph{Inhomogeneous spaces of
positive curvature}, Differential Geom. Appl. \textbf{2} (1992),
123--132.
\bibitem{E2} J.-H. Eschenburg, \emph{Almost positive curvature on the Gromoll-Meyer 7-sphere}, preprint.
\bibitem{GM}D. Gromoll and W.T. Meyer, \emph{An exotic sphere with nonnegative sectional curvature}, Ann. of Math. \textbf{100} (1974), 401-406.
\bibitem{PW} P. Petersen and F. Wilhelm, \emph{Examples of Riemannian manifolds with positive curvature almost everywhere}, Geom. Topol. \textbf{3} (1999), 331-367.
\bibitem{Wallach} N. Wallach, \emph{Compact Riemannian manifolds with strictly positive curvature}, Ann. of Math. \textbf{96} (1972), 277-295.
\bibitem{Wilking} B. Wilking, \emph{Manifolds with positive sectional
curvature almost everywhere}, Invent.\ Math. \textbf{148} (2002), 117--141.
\bibitem{Wil} F. Wilhelm, \emph{An exotic sphere with positive curvature almost everywhere}, Geom. and Top. \textbf{3} (1999), 331-367.
\bibitem{WB} B. Wilking, \emph{The normal homogeneous space $SU(3)\times SO(3)/U(2)$ has positive sectional curvature}, Proceedings of the AMS, \textbf{127} (1999), 1191-1994.
\bibitem{fat} W. Ziller, \emph{Fatness revisited}, lecture notes, University of Pennsylvania, 1999.
\end{thebibliography}

\end{document}